\documentclass{article}

\usepackage{amssymb, latexsym}
\usepackage{graphicx}
\usepackage{amsmath}
\usepackage{color}

\def\reel{{\rm I}\!{\rm R} }

\def\R{\reel}

\def\en t{{{\rm Z}\mkern-5.5mu{\rm Z}}}

\newtheorem{theorem}{Theorem}[section]

\newtheorem{corollary}[theorem]{Corollary}

\newtheorem{lemma}[theorem]{Lemma}

\newtheorem{remark}[theorem]{Remark}

\date{ }

\begin{document}

\title{\Large\bf On the first eigenvalue of the Witten-Laplacian and the diameter of compact shrinking solitons}

\author{Akito Futaki \thanks{Research supported by JSPS Grant-in-Aid for Scientific Research (A) No. 21244003 and Challenging Exploratory Research N0. 23654023. }, \ \
Haizhong Li\thanks{Research supported by NSFC No. 10971110 and Tsinghua University-K. U. Leuven Bilateral Scientific Cooperation Fund.}, \ \
Xiang-Dong Li\thanks{Research supported by NSFC No. 10971032, Key Laboratory RCSDS, CAS, No. 2008DP173182, AMSS Research Grant Y129161ZZ1, and a
Hundred Talents Project of AMSS,
CAS.}
}

\maketitle

\begin{abstract} We prove a lower bound estimate for the first  non-zero eigenvalue of the Witten-Laplacian on  compact
Riemannian manifolds. As an application, we derive a lower bound estimate for the diameter of compact gradient shrinking
Ricci solitons. Our results improve some previous estimates which were obtained by the first author and Y. Sano in \cite{FS},
and by B. Andrews and L. Ni in \cite{AN}. Moreover, we extend the diameter estimate to compact self-similar shrinkers of mean curvature flow.
\end{abstract}

\section{Introduction}

During the recent years, the Bakry-Emery Ricci curvature has received a lot of attention in various areas in mathematics. On the one hand, it has been used to establish some functional inequalities which play an important role in the study of the rate of convergence to the equilibrium measure of diffusion processes \cite{BE, BL1}. On the other hand, it is a good substitute of the Ricci curvature for establishing many interesting theorems in differential geometry, for example, Myers' theorem, eigenvalues estimates, Li-Yau Harnack inequality, Liouville theorems and the Cheeger-Gromoll splitting theorem \cite{BQ1, BQ2, BL2, FLZ, Li05, Q, WW, MW}. Moreover, it has been an important tool in the optimal transport theory \cite{Vil} and in Perelman's work for the entropy formula on Ricci flow \cite{Pel} (see also \cite{Li10}). The purpose of this paper is to prove a lower bound estimate for the first non-zero eigenvalue of the Witten-Laplacian on compact Riemannian manifolds with  Bakry-Emery Ricci curvature bounded from below by a constant. As an application, we prove a lower bound estimate of the diameter for compact shrinking Ricci solitons. Our results improve some previous estimates which were obtained by the first author and Y. Sano in \cite{FS}, and by B. Andrews and L. Ni in \cite{AN}. Moreover, we extend the diameter estimate to compact self-shrinkers of the mean curvature flow.

Let us first introduce some basic notations. Let $(M, g)$ be a complete Riemannian manifold, $f\in C^2(M)$ and $d\mu=e^{-f}dv$, where $dv$ denotes the Riemannian volume measure on $(M, g)$.
For all $u, v\in C^\infty_0(M)$, the following integration by parts formula holds
\begin{eqnarray*}
\int_M \langle \nabla u, \nabla v\rangle d\mu=-\int_M (\Delta_f u)vd\mu=-\int_M u(\Delta_f v)d\mu,
\end{eqnarray*}
where $\Delta_f$ is the so-called Witten-Laplacian on $(M, g)$ with respect to the weighted volume measure $\mu$. More precisely, we have
\begin{eqnarray*}
\Delta_f =\Delta -\nabla f\cdot\nabla. \label{WL}
\end{eqnarray*}

In \cite{BE}, Bakry and Emery proved that for all $u\in C_0^\infty(M)$,
\begin{eqnarray}
\Delta_f |\nabla u|^2-2\langle \nabla u, \nabla \Delta_f u\rangle=2|\nabla^2 u|^2+2(Ric+\nabla^2 f)(\nabla u, \nabla u). \label{BWF}
\end{eqnarray}
The formula $(\ref{BWF})$ can be viewed as a natural extension of the Bochner-Weitzenb\"ock formula. The quantity $Ric+\nabla^2 f$, which is called in the literature the Bakry-Emery Ricci curvature on the weighted Riemannian manifolds $(M, g, f)$, plays as a good substitute of the Ricci curvature in many problems in comparison geometry on weighted Riemannian manifolds. See \cite{BE, BL1, BL2, BQ1, BQ2, FLZ, Li05, Li10, Q, WW, MW} and reference therein.

Now we state the main results of this paper. The first result of this paper is the following lower bound estimate of the first non-zero eigenvalue of the Witten-Laplacian on compact Riemannian manifolds.

\begin{theorem}\label{th1} Let $(M, g)$ be an $n$-dimensional compact Riemannian manifold, and let $\phi\in C^2(M)$. Suppose that there exists a constant $K\in \mathbb{R}$ such that
\begin{eqnarray*}
Ric+\nabla^2\phi\geq Kg.
\end{eqnarray*}
Then the first non-zero eigenvalue $\lambda_1$ of the Witten-Laplacian $\Delta_\phi$ satisfies
\begin{eqnarray}
\lambda_1\geq  \sup\limits_{s\in (0, 1)}\left\{4s(1-s){\pi^2\over d^2}+sK\right\},\label{gap1}
\end{eqnarray}
where $d$ is the diameter of $(M, g)$.
\end{theorem}

As an application of the above theorem, we have the following lower bound estimate for the diameter of compact gradient shrinking Ricci solitons. Recall that a complete Riemannian manifold $(M, g)$ is called a gradient shrinking Ricci soliton if there exists a positive constant $\lambda>0$ and a smooth function $f$ on $M$ such that (see \cite{Cao} and reference therein)
\begin{eqnarray}
Ric(g)+\nabla^2 f=\lambda g. \label{SRS}
\end{eqnarray}
If $f$ is a constant, then $g$ is Einstein. In this case we say that $(M, g, f)$ is trivial.

\begin{theorem}\label{th2} Let $(M, g, f)$ be a non-trivial compact shrinking Ricci soliton with
\begin{eqnarray*}
Ric+\nabla^2 f=\lambda g,
\end{eqnarray*}
where $\lambda$ is a positive constant. Then the diameter of $(M, g)$ satisfies
\begin{eqnarray}
d&\geq& \frac{2(\sqrt{2}-1)\pi} {\sqrt{\lambda}}. \label{d1}
\end{eqnarray}
\end{theorem}

\begin{corollary} Let $(M, g, f)$ be a compact shrinking Ricci soliton with
\begin{eqnarray*}
Ric+\nabla^2 f=\lambda g,
\end{eqnarray*}
where $\lambda$ is a positive constant. If the diameter of $(M, g)$ satisfies
\begin{eqnarray*}
d&<&\frac{2(\sqrt{2}-1)\pi}{\sqrt{\lambda}},
\end{eqnarray*}
then $(M, g)$ must be Einstein.
\end{corollary}

\begin{remark} The study of lower bound estimate of the first eigenvalue on Riemannian manifolds has a long time history. See \cite{Lich, Be, Ch, LY, SY, ZY, CW1, CW2, BQ1} and reference therein. In the case where $(M, g)$ is a compact Riemannian manifold with non-negative Ricci curvature, Zhong and Yang \cite{ZY} obtained the optimal lower bound estimate of the first eigenvalue of the Laplacian, i.e.,
\begin{eqnarray*}
\lambda_1\geq {\pi^2\over d^2},
\end{eqnarray*}
where $d$ is the diameter of $(M,g)$. In \cite{SZ07},  Shi and Zhang proved that on compact Riemannian manifolds with Ricci curvature bounded below by a constant $K\in \mathbb{R}$, i.e.,
\begin{eqnarray*}
Ric\geq Kg.
\end{eqnarray*}
the first non-zero eigenvalue $\lambda_1$ of the Laplacian $\Delta$ satisfies
\begin{eqnarray*}
\lambda_1\geq  \sup\limits_{s\in (0, 1)}\left\{4s(1-s){\pi^2\over d^2}+sK\right\}.
\end{eqnarray*}
See also Qian-Zhang-Zhu \cite{QZZ} for its extension to compact Alexandrov spaces with Ricci curvature bounded from below by
$K$. Note that, by direct calculation, we have
\begin{equation*}
\sup\limits_{s\in (0, 1)}\left\{4s(1-s){\pi^2\over d^2}+sK\right\}
=\left\{\begin{array}{ll}
0 &{\rm if}\ \ Kd^2<-4\pi^2,\\
\left({\pi\over d}+{Kd\over 4\pi}\right)^2&{\rm if}\ \ Kd^2\in [-4\pi^2, 4\pi^2],\\
K &{\rm if}\ \ Kd^2\in (4\pi^2, (n-1)\pi^2].
\end{array}\right.
\end{equation*}
Theorem \ref{th1} is a natural extension of the above estimate of Shi and Zhang to the Witten-Laplacian via the Bakry-Emery Ricci curvature on compact Riemannian manifolds, and it is a general principle that
this kind of extension is possible. A typical such result can also be found in \cite{LuRowlettev3}.
\end{remark}

\begin{remark} In \cite{FS}, under the same condition as in Theorem \ref{th1}, the first author and Y. Sano proved that the first non-zero eigenvalue of the Witten-Laplacian $\Delta_\phi$ satisfies
\begin{eqnarray}
\lambda_1\geq {\pi^2\over d^2}+{31\over 100}K. \label{gap2}
\end{eqnarray}
One can easily check that our new lower bound estimate $(\ref{gap1})$ for $\lambda_1$ is better than $(\ref{gap2})$. In the case where $(M, g, f)$ is a non-trivial compact gradient shrinking Ricci solitons with $Ric+\nabla^2 f=\lambda g$,  the estimate $(\ref{gap2})$ led to the following lower bound estimate for the diameter of $(M, g)$:
\begin{eqnarray}
d\geq {10\pi\over 13\sqrt{\lambda}}.\label{d2}
\end{eqnarray}
One can easily check that
\begin{eqnarray*}
2(\sqrt{2}-1)&>&\frac{10}{13}.
\end{eqnarray*}
Hence, the lower bound  estimate $(\ref{d1})$ in Theorem \ref{th2} is sharper than $(\ref{d2})$.
\end{remark}

\begin{remark} Taking $s={1\over 2}$ in Theorem \ref{th1}, we obtain the following lower bound estimate
\begin{eqnarray}
\lambda_1\geq {\pi^2\over d^2}+{K\over 2}.\label{s=1/2}\label{gap3}
\end{eqnarray}
This recaptures the lower bound estimate of $\lambda_1$ due to Andrews and Ni (Proposition 3.1 in  \cite{AN}).  In the case where $(M, g, f)$ is a non-trivial compact gradient shrinking Ricci solitons with $Ric+\nabla^2 f=\lambda g$,  Andrews and Ni (Corollary 3.1 in \cite{AN}) used the estimate $(\ref{gap3})$ to derive the following lower bound estimate for the diameter of $(M, g)$:
\begin{eqnarray}
d\geq \sqrt{2\over 3\lambda}\pi,\label{d3}
\end{eqnarray}
which is better than $(\ref{d2})$ obtained in \cite{FS}.
One can easily check that
\begin{eqnarray*}
2(\sqrt{2}-1)&>&\sqrt{\frac 2 3}.
\end{eqnarray*}
Hence, the lower bound  estimate $(\ref{d1})$ in Theorem \ref{th2}  is sharper than $(\ref{d3})$.
\end{remark}

As another application of Theorem \ref{th1}, we can also obtain a lower bound estimate for the diameter of compact self-shrinkers
of the mean curvature flow (Theorem \ref{th4}). In fact Theorem \ref{th2} above follows from Theorem \ref{th1} and the fact that
the Witten-Laplacian on the Ricci soliton with (\ref{SRS}) has eigenvalue $2\lambda$ (\cite{FS}, see Lemma \ref{lfs} below), but
we can show that the Witten-Laplacian takes the same eigenvalue $2\lambda$ on a compact self-shrinker
\begin{equation}\label{mcfl}
x^\perp = - \frac 1 \lambda \vec{H}.
\end{equation}

\section{Proof of Theorem \ref{th1}}

To prove Theorem \ref{th1}, we use the following comparison theorem due to Chen and Wang \cite{CW1, CW2},  Bakry and Qian \cite{BQ1}, also Andrews and Ni \cite{AN}.

\begin{theorem}\label{th4} (Chen-Wang \cite{CW1, CW2}, Bakry-Qian \cite{BQ1}, Andrews-Ni \cite{AN}) Let $(M, g)$ be an $n$-dimensional compact Riemannian manifold, and let $\phi\in C^2(M)$. Suppose that there exists a constant $K\in \mathbb{R}$ such that
\begin{eqnarray*}
Ric+\nabla^2\phi\geq Kg.
\end{eqnarray*}
Then the first non-zero Neumann eigenvalue of the Witten-Laplacian $\Delta_\phi$ satisfies
\begin{eqnarray}
\lambda_1\geq \lambda_1(L),\label{com1}
\end{eqnarray}
where $\lambda_1(L)$ denotes the first non-zero Neumann eigenvalue of the following one-dimensional Ornstein-Uhlenbeck operator on $ (-{d/2}, {d/2})$:
\begin{eqnarray}
L={d^2\over dx^2}-Kx{d\over dx}.\label{OUL}
\end{eqnarray}
More precisely, $\lambda_1(L)$ is the first non-zero eigenvalue of the problem
\begin{eqnarray*}
& &v''(x)-Kxv'(x)=-\lambda v(x), \ \ \ \ x\in (-{d\over 2}, {d\over 2}),\\
& &\ \ \ \ \ \ \ \ \ \ v'(-{d\over 2})=0, \ \  \ \ v'({d\over 2})=0.
\end{eqnarray*}
\end{theorem}

\noindent{\bf Proof of Theorem \ref{th1}}.  By Theorem \ref{th4}, we need only to prove that
\begin{eqnarray}
\lambda_1(L)\geq \sup\limits_{s\in (0, 1)}\left\{4s(1-s){\pi^2\over d^2}+sK\right\}.\label{aaa}
\end{eqnarray}
To prove $(\ref{aaa})$, we modify the argument used in the proof of Theorem 1.1  in \cite{SZ07}, cf. also the one of Corollary 4.3 in \cite{QZZ}. Denote $D={d\over 2}$, $f=v'$. Then $f$ is the eigenfunction of the first non-zero eigenvalue $\lambda-K$ of the Ornstein-Uhlenbeck operator $L$ on $(-D, D)$ with the Dirichlet boundary condition. More precisely, we have
\begin{eqnarray}
f''-Kx f'=-(\lambda-K) f,\label{D}.
\end{eqnarray}
and $f(-D)=f(D)=0$.
By maximum principle, we can prove that  $f$ must have fixed sign on $(-D, D)$. So we can assume that $f(x)>0$ for all $x\in (-D, D)$.

Fix a constant $a>1$. Multiplying $f^{a-1}$ to the both sides of $(\ref{D})$ and integrating on $(-D, D)$, we have
\begin{eqnarray*}
\int_{-D}^D f^{a-1}(x)f''(x)dx=-(\lambda-K)\int_{-D}^D f^a(x)dx+\int_{-D}^D Kxf^{a-1}(x)f'(x)dx.
\end{eqnarray*}
Integrating by parts and using the fact $f(\pm D)=0$ and $f(x)>0$ on $(-D, D)$, we get
\begin{eqnarray*}
\int_{-D}^D f^{a-1}(x)f''(x)dx=-(a-1)\int_{-D}^D f^{a-2}(x)f'^2(x)dx=-{4(a-1)\over a^2}\int_{-D}^D [(f^{a/2}(x))']^2dx,
\end{eqnarray*}
Similarly, we have
\begin{eqnarray*}
\int_{-D}^D  Kxf^{a-1}(x)f'(x)dx&=&-\int_{-D}^D f(Kf^{a-1}+(a-1)Kxf^{a-2}f')dx\\
&=&-K\int_{-D}^D f^{a}(x)dx-(a-1)\int_{-D}^D Kxf^{a-1}(x)f'(x)dx,
\end{eqnarray*} which yields
\begin{eqnarray*}
\int_{-D}^D  Kx f^{a-1}(x)f'(x)dx&=&-{K\over a}\int_{-D}^D f^{a}(x)dx.
\end{eqnarray*}
Set $u=f^{a/2}$. Then we can deduce that
\begin{eqnarray*}
{4(a-1)\over a^2}\int_{-D}^D |u'|^2dx%&=&(\lambda-K)\int_{-D}^D u^2dx+{K\over a}\int_{-D}^D u^2dx\\
=\left(\lambda-K(1-{1\over a})\right)\int_{-D}^D u^2dx.
\end{eqnarray*}
Let $s=1-{1\over a}$. We get
\begin{eqnarray*}
4s(1-s)\int_{-D}^D |u'|^2dx=(\lambda-Ks)\int_{-D}^D u^2dx.
\end{eqnarray*}
Hence
\begin{eqnarray*}
{\lambda-Ks\over 4s(1-s)}={\int_{-D}^D |u'|^2dx\over \int_{-D}^D u^2dx}.
\end{eqnarray*}
Using the Wirtinger inequality, we have
\begin{eqnarray*}
{\lambda-Ks\over 4s(1-s)}\geq {\pi^2 \over 4D^2}={\pi^2\over d^2}.
\end{eqnarray*}
Hence, for all $s\in (0, 1)$, we have
\begin{eqnarray*}
\lambda\geq 4s(1-s){\pi^2\over d^2}+Ks.
\end{eqnarray*}
The proof of Theorem \ref{th1} is completed. \hfill $\square$

\section{Proof of Theorem \ref{th2}}

To prove Theorem \ref{th2}, we need the following

\begin{lemma}\label{lfs}(\cite{FS}) Let $(M, g, f)$ be a non-trivial gradient shrinking Ricci solitons with
\begin{eqnarray}
Ric+\nabla^2 f=\lambda g.
\end{eqnarray}
Then $f$ is an eigenfunction of the Witten-Laplacian $\Delta_f$ with eigenvalue equal to $2\lambda$.
\end{lemma}
\medskip

\noindent{\it Proof of Theorem \ref{th2}}.  By Theorem \ref{th1} and Lemma \ref{lfs}, for all $s\in (0, 1)$, we have
\begin{eqnarray*}
2\lambda\geq 4s(1-s){\pi^2\over d^2}+s\lambda,
\end{eqnarray*}
which yields
\begin{eqnarray*}
\lambda \ge {4s(1-s) \over 2-s} \frac{\pi^2}{d^2}.
\end{eqnarray*}
But elementary computations show
$$ \frac{4s(1-s)}{2-s} \le 12 - 8\sqrt{2},$$
where the equality is attained for $s = 2 - \sqrt{2} \in (0,1)$.  Thus we obtain
%Let
%\begin{eqnarray*}
%f(s):={4\pi^2\over d^2} s^2-\left(\lambda+{4\pi^2\over d^2}\right)s+2\lambda.
%\end{eqnarray*}
%Then $f(s)\geq 0$ holds for all $s\in (0, 1)$. Note that
%\begin{eqnarray*}
%f(s)&=&{4\pi^2\over d^2}\left(s-{1\over 2}-{d^2\lambda\over 8\pi^2}\right)^2+2\lambda-{\pi^2\over d^2}\left(1+{d^2\lambda\over 4\pi^2}\right)^2.
%\end{eqnarray*}
%By direct calculation, we have
%\begin{eqnarray*}
%\inf\limits_{s\in (0, 1)}f(s)= f\left({1\over 2}+{d^2\lambda\over 8\pi^2}\right)={3\lambda\over 2}-{\pi^2\over d^2}-{d^2\lambda^2\over 4^2\pi^2}.
%\end{eqnarray*}
%Therefore $\lambda$ must satisfy the constraint condition
%\begin{eqnarray*}
%{3\lambda\over 2}-{\pi^2\over d^2}-{d^2\lambda^2\over 4^2\pi^2}\geq 0.
%\end{eqnarray*}
%This yields
\begin{eqnarray*}
\lambda\geq 4(3-2\sqrt{2}) {\pi^2\over d^2}.
\end{eqnarray*}
Equivalently, the diameter of $(M, g)$ must satisfy the following lower bound
\begin{eqnarray*}
d&\geq& \frac{2(\sqrt{2}-1)\pi}{\sqrt{\lambda}}.
\end{eqnarray*}
The proof of Theorem \ref{th2} is completed. \hfill $\square$

\section{The diameter of compact self-shrinkers for mean curvature flow}

Let $x:M \to \R^{n+p}$ be an $n$-dimensional submanifold in
the (n+p)-dimensional Euclidean space. If we let the position
vector $x$ evolve in the direction of the mean curvature $\vec{H}$,
then it gives rise to a solution to the mean curvature flow :
$$
x:M\times [0,T)\rightarrow \R^{n+p}, \qquad \frac{\partial
x}{\partial t} = \vec{H}.
$$

We call the immersed manifold $M$ a self-shrinker if it satisfies
the quasilinear elliptic system (see \cite{H}, or \cite{CM}): for some positive constant $\lambda$,
$$
\vec{H}=-\lambda x^{\perp},
$$
where $\perp$ denotes the projection onto the normal bundle of $M$.

We have (see \cite{LW})
$$
\frac 1 {2\lambda} |\vec{H}|^2+\frac{1}{4}\Delta |x|^2=\frac{n}{2}.
$$
Put
$$\phi := 2\lambda(\frac{|x|^2}4 -\frac{n}{4\lambda}).
\eqno (15)
$$
Define the Witten-Laplacian by
$$
\Delta_\phi=\Delta- \nabla\phi\cdot\nabla.
$$
From above formulas, we can check
$$
\begin{array}{lcl}
\Delta_\phi(\frac{1}{4}|x|^2)&=&\Delta (\frac{1}{4}|x|^2)-\frac{\lambda}{8}\nabla|x|^2\cdot \nabla|x|^2\\
&=&\frac{n}{2}- \frac 1{2\lambda}|\vec{H}|^2-\frac{\lambda}2|x^T|^2\\
&=&\frac{n}{2}-\frac{\lambda}{2}|x|^2.
\end{array}
$$
\noindent
Thus we have
$$
\Delta_\phi(\frac{1}{4}|x|^2-\frac{n}{4\lambda})=
-2\lambda(\frac{|x|^2}4-\frac{n}{4\lambda}).
$$

Thus we have proved

\begin{theorem}\label{th3} In the above situation we have the eigenvalue $2\lambda$ of  the Witten-Laplacian $\Delta_\phi$ with
 eigenfunction $\phi$:
$$
\Delta_\phi\phi=
-2\lambda \phi.
$$
\end{theorem}

Let $h^\alpha_{ij}$ is the components of the second fundamental form, $H^\alpha=\sum\limits_kh^\alpha_{kk}$, the Gauss equation is (see \cite{CL})

$$
R_{ij}=\sum\limits_\alpha H^\alpha h^\alpha_{ij}-\sum\limits_{\alpha,k}h^\alpha_{ik}h^\alpha_{kj}.
$$
Thus we have from the definition of $\phi$ in (15)
$$
\begin{array}{lcl}
R_{ij}+\phi_{ij}&=&\lambda g_{ij}-\sum\limits_{\alpha,k}h^\alpha_{ik}h^\alpha_{kj}\\
&\geq& [\lambda -K_0]g_{ij},
\end{array}
$$
where
$$
K_0=\max_{1\leq i\leq
n}[\sum\limits_{\alpha,k}h^\alpha_{ik}h^\alpha_{ki}],
$$
and we have used
$$
\phi_{ij}=(\frac{\lambda}{2}|x|^2)_{ij}=\lambda g_{ij}-\sum\limits_\alpha H^\alpha h^\alpha_{ij}.
$$
By the similar argument as previous sections, we have
\begin{theorem}\label{th3} Let $X:M\to R^{n+p}$ be an $n$-dimensional compact self-shrinker.
Suppose that there exists a constant $K\in \mathbb{R}$ such that
$$
Ric+\nabla^2\phi\geq Kg,
$$
where
$$
K=\lambda -K_0,\qquad K_0= \max_{1\leq i\leq n}[\sum\limits_{\alpha,k}h^\alpha_{ik}h^\alpha_{ki}].
$$
Then the first non-zero eigenvalue $\lambda_1$ of the Witten-Laplacian $\Delta_\phi$ satisfies
$$
\lambda_1\geq  \sup\limits_{s\in (0, 1)}\left\{4s(1-s){\pi^2\over d^2}+sK\right\},
$$
where $d$ is the diameter of $M$.
\end{theorem}

Following the arguments of previous sections, we have
$$
2\lambda \geq
4s(1-s)\frac{\pi^2}{d^2}+sK
\eqno (16)
$$
for all $s\in (0,1)$.

Thus we obtain a diameter estimate for compact self-shrinker, which
are not minimal submanifold of $S^{n+p-1}(\sqrt{n/\lambda})$ (which corresponds $|x|^2=constant$,
so from (15) we get  $\phi=0$, trivial case).
Choosing $s=\frac{1}{2}$ in (16), we have for $n$-dimensional self-shrinkers

\begin{theorem}\label{th4} Let $x:M\to R^{n+p}$ be an $n$-dimensional compact self-shrinker such that $x(M)$ is not minimal
submanifold in  $S^{n+p-1}(\sqrt{n/\lambda})$, and let
$h^\alpha_{ij}$ be the components of the the second fundamental form
of $M$. Then we have
$$
d\geq \frac{1}{\sqrt{\frac{3\lambda}{2}+\frac{1}{2}K_0}}\pi,
$$
where
$$
K_0:=\max_{1\leq i\leq n}[\sum\limits_{\alpha,k}h^\alpha_{ik}h^\alpha_{ki}].
$$
\end{theorem}

When $p=1$, we have
\begin{corollary}\label{cor} Let $x:M\to R^{n+1}$ be an $n$-dimensional compact self-shrinker such that $x(M)$ is not
$S^n(\sqrt{n/\lambda})$,  and let $\lambda_i$ be the principal
curvatures of $M$.  Then we have
$$
d\geq \frac{1}{\sqrt{\frac{3\lambda}{2}+\frac{1}{2}K_0}}\pi,
$$
where
$$
K_0:= \max_{p \in M}\max_{1\leq i\leq n}\lambda_i^2.
$$
\end{corollary}

\begin{remark} The similar results hold for self-shrinkers in Riemannian cone manifolds as considered in \cite{FHY12}.
\end{remark}

\bigskip

\noindent{\bf Acknowledgement}. Part of this work was done during the first author visited Mathematical Sciences Center of Tsinghua University in September-October 2011, by the invitation of Professor S.-T. Yau. The third author would like to thank Professors Dominique Bakry and Mu-Fa Chen for explaining their results obtained in \cite{BQ1, CW1, CW2} to him. We also thank Dr. Daguang Chen for helpful discussion.

\medskip

\begin{flushleft}
\medskip\noindent
Akito Futaki, {\sc Department of Mathematics, Tokyo Institute of Technology,
O-okayama, Meguro, Tokyo, 152-8551, Japan,} \\
E-mail: futaki@math.titech.ac.jp

\medskip

Haizhong Li, {\sc Department of Mathematical Sciences, Tsinghua University, Beijing, 100084, China,}\\
E-mail: hli@math.tsinghua.edu.cn

\medskip

Xiang-Dong Li, {\sc  Academy of Mathematics and System Science, Chinese
Academy of Sciences, 55, Zhongguancun East Road, Beijing, 100190, P. R. China,}\\
E-mail: xdli@amt.ac.cn
\end{flushleft}

\end{document}